\documentclass{article}
\usepackage[utf8]{inputenc}
\usepackage{float}
\usepackage{amsmath}
\usepackage{amsthm}
\usepackage{amssymb}
\usepackage{graphicx}
\usepackage{setspace}
\usepackage[bookmarks=true,bookmarksnumbered=true,bookmarksopen=true,bookmarksopenlevel=1,
 breaklinks=false,pdfborder={0 0 1},backref=false,colorlinks=false]
 {hyperref}

\makeatletter
\numberwithin{equation}{section}
\numberwithin{figure}{section}
\theoremstyle{plain}
\newtheorem{thm}{\protect\theoremname}[section]
\theoremstyle{remark}
\newtheorem{rem}[thm]{\protect\remarkname}
\theoremstyle{definition}
\newtheorem{defn}[thm]{\protect\definitionname}
\theoremstyle{plain}
\newtheorem{lem}[thm]{\protect\lemmaname}

\@ifundefined{date}{}{\date{}}

\usepackage{graphicx}
\usepackage{amsfonts}
\usepackage{lmodern}
\usepackage{amsthm}
\allowdisplaybreaks
\usepackage{cases}

\def\XXint#1#2#3{{\setbox0=\hbox{$#1{#2#3}{\int}$ }
\vcenter{\hbox{$#2#3$ }}\kern-.6\wd0}}

\usepackage{accents}

\usepackage{extarrows}

\makeatother

\usepackage[bibstyle=alphabetic,citestyle=alphabetic-verb]{biblatex}
\providecommand{\definitionname}{Definition}
\providecommand{\lemmaname}{Lemma}
\providecommand{\remarkname}{Remark}
\providecommand{\theoremname}{Theorem}

\AtBeginDocument{ 
\DeclareBibliographyAlias{article}{std}
\DeclareBibliographyAlias{book}{std}
\DeclareBibliographyAlias{booklet}{std}
\DeclareBibliographyAlias{collection}{std}
\DeclareBibliographyAlias{inbook}{std}
\DeclareBibliographyAlias{incollection}{std}
\DeclareBibliographyAlias{inproceedings}{std}
\DeclareBibliographyAlias{manual}{std}
\DeclareBibliographyAlias{misc}{std}
\DeclareBibliographyAlias{online}{std}
\DeclareBibliographyAlias{patent}{std}
\DeclareBibliographyAlias{periodical}{std}
\DeclareBibliographyAlias{proceedings}{std}
\DeclareBibliographyAlias{report}{std}
\DeclareBibliographyAlias{thesis}{std}
\DeclareBibliographyAlias{unpublished}{std}
\DeclareBibliographyAlias{*}{std}

\DeclareBibliographyDriver{std}{%
  \usebibmacro{bibindex}%
  \usebibmacro{begentry}%
  \usebibmacro{author/editor+others/translator+others}%
  \setunit{\labelnamepunct}\newblock
  \usebibmacro{title}%
  \newunit\newblock
  \printfield{journaltitle}
  \newunit\newblock
  \printfield{publisher}
  \newunit\newblock
  \printfield{isbn}
  \newunit\newblock
  \printfield{url}
  \newunit\newblock
  \usebibmacro{date}%
  \newunit\newblock
  \usebibmacro{finentry}}
}

\begin{document}
\title{A simple geometric construction of an ODE with undecidable blow-ups}
\author{Manh Khang Huynh}
\date{\today}
\maketitle
\begin{abstract}
We present a simple construction of an ODE on $\mathbb{R}^{n}$ where
the vector field is smooth, and finite-time blow-up is equivalent
to the halting problem for a universal Turing machine.
\end{abstract}

\section{Introduction}

The classical ODE problem 
\begin{align}
\dot{x}\left(t\right) & =F\left(x\left(t\right)\right)\label{eq:ODE}\\
x\left(0\right) & =x_{0}\nonumber 
\end{align}
where $F:\mathbb{R}^{n}\to\mathbb{R}^{n}$ is smooth, admits a unique
smooth solution $x\left(t\right)$ which lives on a maximal time interval
of existence $\left(t_{\min},t_{\max}\right)$. In particular if $t_{\max}<\infty$
then we must have finite-time blow-up \cite{teschlOrdinaryDifferentialEquations2012},
i.e. $\lim_{t\uparrow t_{\max}}\left|x\left(t\right)\right|=$$\infty$.
Given $F$ and $x_{0}$, determining the value of $t_{\max}$, or
simply whether it is finite, is an unsolved problem. In certain scenarios,
we can use numerical methods to rigorously verify that the ODE solution
will blow up in some specific direction \cite{takayasuNumericalValidationBlowsolutions2017}. 

In \cite{taoFiniteTimeBlowup2016a}, Terence Tao constructed a solution
to an averaged version of 3D Navier-Stokes which exhibits finite-time
blow-up. The ideas behind the ``logic gates'' used in the construction
led to \cite{taoUniversalityPotentialWell2017}, in which the author
used Nash embedding to embed any vector-field flows with an ``adapted
1-form'', including a smooth model of a universal Turing machine,
into a potential well system (which is an ODE on the cotangent bundle
of $\mathbb{R}^{n}$). In particular, it also implied the halting
problem was equivalent to the trajectory entering some bounded open
set, with the initial position corresponding to the input tape for
the Turing machine (see also \parencite{mooreUnpredictabilityUndecidabilityDynamical1990,mooreGeneralizedShiftsUnpredictability1991}
for the theory of generalized shifts). Similar Turing completeness
results were then proved for steady Euler flows, and we refer to \parencite{cardonaUniversalityEulerFlows2019,cardonaConstructingTuringComplete2020,cardonaTuringUniversalityIncompressible2022}
for more details.

Terence Tao speculated in \cite{taoUniversalityPotentialWell2017}
that by Turing completeness, one could possibly ``program'' a solution
towards blow-up. In this paper we present such an ODE, where the aforementioned
open set becomes ``infinity''.
\begin{thm}[Main theorem]
\label{thm:main}  We can find $n\in\mathbb{N}$ and a smooth function
$F:\mathbb{R}^{n}\to\mathbb{R}^{n}$, with an associated universal
Turing machine $\mathrm{TM}$, and a set of initial positions $X_{0}$,
such that for any $x_{0}\in X_{0}$ (corresponding to some input tape
$s_{0}$), the solution in \eqref{eq:ODE} exhibits (finite-time)
blow-up \footnote{an ODE solution could ``grow up'', i.e. tend to infinity at infinite
time, for instance in \cite{takayasuNumericalValidationBlowsolutions2017},
but such a scenario is ruled out in our construction.} if and only if the Turing machine $\mathrm{TM}$ halts when given
the input tape $s_{0}$.
\end{thm}

\begin{rem}
The undecidability of ODE blow-up seems natural and expected, given
the possible complexity of ODE dynamics. In particular, it was proven
in \cite{gracaBoundednessDomainDefinition2008,gracaComputabilityNoncomputabilityUndecidability2009}
for analytic vector fields, using computability theory techniques to robustly simulate a Turing machine with small error bounds. Our
construction was independently obtained via a short and simple geometric embedding scheme, with an exact simulation of a Turing machine. 
\end{rem}

\subsubsection*{Acknowledgements}

The author would like to thank Jean-Philippe Lessard and Alex Blumenthal
for interesting discussions about the theory of rigorous numerical
verification.

\section{Preliminaries}

We first recall the definition of a Turing machine.
\begin{defn}
A Turing machine $\left(Q,\mathrm{START},\mathrm{HALT},\Sigma,\delta\right)$
consists of
\begin{itemize}
\item A finite set of states $Q$, where $\{\mathrm{START},\mathrm{HALT}\}\subseteq Q$.
\item A finite alphabet set $\Sigma$ with at least 2 elements.
\item A transition function $\delta:Q\backslash\left\{ \mathrm{HALT}\right\} \times\Sigma\to Q\times\Sigma\times\left\{ -1,0,1\right\} $
\end{itemize}
We denote $q$ as the current state and $t=\left(t_{n}\right)_{n\in\mathbb{Z}}\in\Sigma^{\mathbb{Z}}$
as the current tape. The machine works by following a simple algorithm:
\begin{enumerate}
\item Set $q$ to $\mathrm{START}$ and $t$ to the input tape.
\item If $q$ is $\mathrm{HALT}$, halt and return $t$ as output. Otherwise
compute $\left(q',t'_{0},\epsilon\right)=\delta\left(q,t_{0}\right)$.
\item Replace $q$ with $q'$ and $t_{0}$ with $t'_{0}$.
\item Replace the tape with its $\epsilon-$shifted version (with $\epsilon=1$
being left-shift, $0$ being no shift, and $-1$ being right-shift).
Then return to step 2.
\end{enumerate}
Then for any $q\in Q\backslash\left\{ \mathrm{HALT}\right\} $ and
$t\in\Sigma^{\mathbb{Z}}$, we have an auxiliary map $\Delta:\left(q,t\right)\mapsto\left(q',t'\right)$
giving the next state $q'$ and the next input tape $t'$. 

From this point on, we fix $\mathrm{TM}$ to be a \emph{universal}
Turing machine \parencite{aroraComputationalComplexityModern2016},
which can be used to model any other Turing machine. 

We recall the standard torus model for $\mathrm{TM}$ as constructed
in \parencite[Section 4]{taoUniversalityPotentialWell2017}. 
\end{defn}

\begin{lem}
Let $M=\mathbb{T}^{2}\times\mathbb{T}^{2}$ be the 4-dimensional torus.
Each state $q\in Q$ corresponds to an open square $B_{q}\subseteq\mathbb{T}^{2}$
(with the closed squares $\left(\overline{B_{q}}\right)_{q\in Q}$
being disjoint) , while each string $s\in\Sigma^{\mathbb{Z}}$ corresponds
to a unique point $f_{s}\in\mathbb{T}^{2}$, and there is a diffeomorphism
$\Phi:M\to M$ such that $\Phi\left(\overline{B_{q}}\times\left\{ f_{s}\right\} \right)\subseteq B_{q'}\times\left\{ f_{s'}\right\} $
for any $\left(q,t\right)$ where $\Delta\left(q,t\right)=\left(q',t'\right)$.
\end{lem}

\subsubsection*{Halting region $U_{0}$}

With the torus model constructed above, we now define $x_{\mathrm{START}}$
as the center of $B_{\mathrm{START}}$, and the open set $U_{0}\subseteq M$
as a small neighborhood of $\overline{B_{\mathrm{HALT}}}\times\mathbb{T}^{2}$,
such that $\overline{U}_{0}$ is disjoint from any $\overline{B_{q}}\times\mathbb{T}^{2}$
where $q\neq\mathrm{HALT}$. Then given a starting point $\left(x_{\mathrm{START}},f_{s}\right)$,
the trajectory $\Phi^{\mathbb{N}}\left(x_{\mathrm{START}},f_{s}\right)$
enters $\overline{B_{\mathrm{HALT}}}\times\mathbb{T}^{2}$ at some
time if and only if the Turing machine $\mathrm{TM}$ halts when given
input tape $s$. A similar statement holds for the trajectory entering
$U_{0}$. We could also modify $U_{0}$ to specify some output tape
digits. $U_{0}$ can also be chosen to be connected with smooth boundary,
and we can either shrink or grow $U_{0}$ without changing the conclusion.

\subsubsection*{Continuous-time model}

To turn discrete-time dynamics into continuous-time dynamics, we consider
$M_{I}$ which is $M\times\left[0,1\right]/\sim$ where each $\left(p,1\right)$
is identified with $\left(\Phi\left(p\right),0\right)$. We parametrize
$M_{I}$ by variables $\left(p,t\right)$, and consider the one-parameter
group, or flow, $\left(\Phi_{Y}^{t}\right)_{t\in\mathbb{R}}$ generated
by the vector field $Y=\partial_{t}$. Then there is no stationary
point (equilibrium). Now we let $U_{1}$ be a small neighborhood containing
$U_{0}\times\left\{ 0\right\} $ (disjoint from any $\left(\overline{B_{q}}\times\mathbb{T}^{2}\right)\times\left[0,1\right]$
for $q\neq\mathrm{HALT}$). Obviously $U_{1}$ can be chosen to be
connected with smooth boundary. 

Consequently, given the starting point $\left(x_{\mathrm{START}},f_{s},0\right)\in M_{I}$,
the trajectory along $Y$ enters $U_{1}$ at some time if and only
if the Turing machine $\mathrm{TM}$ halts when given the input tape
$s$.

\section{Geometric embedding}

\begin{figure}[H]
\centering
\includegraphics[width=0.4\textwidth]{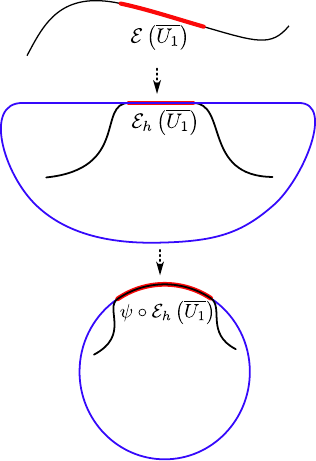}

\caption{Geometric embedding scheme}
\end{figure}

We can use Whitney embedding (or Nash embedding) to smoothly embed
the manifold $M_{I}$ constructed above into a compact region of $\mathbb{R}^{d}$
(where $d$ could be chosen to be $10$, for instance). Let $\mathcal{E}$
be this embedding. 

\subsubsection*{Height function}

There is a smooth function $h$ on $M_{I}$ such that $0\leq h\leq1$
and $\overline{U_{1}}=h^{-1}\left(\left\{ 1\right\} \right)$. We
can therefore define the embedding 
\begin{align*}
\mathcal{E}_{h}:M_{I} & \to\mathbb{R}^{d+1}\\
p & \mapsto\left(\mathcal{E}\left(p\right),h\left(p\right)\right)
\end{align*}
Let $M_{h}$ be the image of $\mathcal{E}_{h}$. Then $\mathcal{E}_{h}\left(\overline{U_{1}}\right)=\left\{ x_{d+1}=1\right\} \cap M_{h}$.
The pushforward $\left(\mathcal{E}_{h}\right)_{*}Y$ of $Y$ can be
smoothly extended to a vector field $Y_{h}$ on $\mathbb{R}^{d+1}$
with support near $M_{h}$. Consequently, reaching $\left\{ x_{d+1}=1\right\} $
via the flows of $Y_{h}$ is equivalent to the halting problem.

\subsubsection*{Spherical embedding}

As $M_{h}$ is compact and bounded, there is a compact manifold $B$
with smooth boundary such that 
\begin{enumerate}
\item $M_{h}\subseteq B\subseteq\mathbb{R}^{d+1}$.
\item $M_{h}\cap\partial B=\left\{ x_{d+1}=1\right\} \cap M_{h}=\mathcal{E}_{h}\left(\overline{U_{1}}\right)$.
\item There is a diffeomorphism $\psi:\mathbb{R}^{d+1}\to\mathbb{R}^{d+1}$
such that $\psi\left(B\right)=\overline{\mathbb{D}^{d+1}}$ (the closed
unit ball of $\mathbb{R}^{d+1}$).
\end{enumerate}
Then
\[
\psi\circ\mathcal{E}_{h}\left(\overline{U_{1}}\right)=\mathbb{S}^{d}\cap\psi\left(M_{h}\right),
\]
so reaching the sphere via the flows of $\psi_{*}Y_{h}$ is equivalent
to the halting problem.

\subsubsection*{Poincaré radial compactification}

Finally, there is a diffeomorphism 
\begin{align*}
T:\mathbb{D}^{d+1} & \to\mathbb{R}^{d+1}\\
x & \mapsto\frac{x}{\sqrt{1-\left|x\right|^{2}}}
\end{align*}
 between the open unit ball and the whole space, which leads to the
Poincaré radial compactification of $\mathbb{R}^{d+1}$, by essentially
``mapping'' the sphere $\mathbb{S}^{d}$ to $\partial\mathbb{R}^{d+1}$. 

We are done. Given the starting point $T\circ\psi\circ\mathcal{E}_{h}\left(x_{\mathrm{START}},f_{s},0\right)$,
the trajectory via the flow of $\left(T\circ\psi\right)_{*}Y_{h}$
blows up in finite time if and only if the Turing machine $\mathrm{TM}$
halts when given the input tape $s$. Thus Theorem \ref{thm:main}
is proven. \qed
\begin{rem}
Since a universal Turing machine could model a halting Turing machine,
there exist initial conditions that lead to finite-time blow-ups.
Similarly, there are initial conditions that would not lead to blow-ups.
We have essentially geometrically ``programmed'' a solution towards
blow-up.
\end{rem}

\begin{rem}
In \parencite{takayasuNumericalValidationBlowsolutions2017}, Poincaré
radial compactification was used to regularize ODE blow-ups for rigorous
numerical verification. See also \parencite{lessardGeometricCharacterizationUnstable2021,lessardSaddleTypeBlowUpSolutions2022}
for similar geometric techniques related to blow-ups. Our construction is flexible and could be further simplified by condensing multiple steps into a single embedding, but the connection to numerical verification and Poincaré radial compactification deserves to be made explicit. 
\end{rem}

\section{Open questions}

Following \cite{taoUniversalityPotentialWell2017}, it is natural
to ask whether we can find a potential well system or a nonlinear
wave system where blow-up is undecidable, using geometric methods.

In yet another direction, we look towards other PDEs such as those
in fluid mechanics. Instead of encoding the Turing-complete dynamical
complexity into the nonlinearity $F$ as in \eqref{eq:ODE}, we can
encode the complexity into the initial velocity field, such as the
steady Euler flows constructed in \parencite{cardonaConstructingTuringComplete2020,cardonaUniversalityEulerFlows2019,cardonaTuringUniversalityIncompressible2022}.
See also \parencite{taoUniversalityIncompressibleEuler2018} for a
discussion of the universality of non-steady Euler flows on manifolds.

\printbibliography

\end{document}